\documentclass{article}
\usepackage{amscd}
\usepackage{amsmath}
\usepackage{amssymb}
\usepackage{amsthm}
\usepackage{ascmac}
\usepackage{graphicx}
\usepackage{here}
\usepackage{mathtools}\numberwithin{equation}{section}
\usepackage{minitoc}
\usepackage{multicol}
\usepackage{multirow}
\usepackage{setspace}
\usepackage[dvipdfmx]{color}
\theoremstyle{definition}
\newtheorem{df}{Definition}
\newtheorem{note}{Notation}
\newtheorem{q}{Problem}
\newtheorem{eg}{Example}
\newtheorem{rmk}{Remark}
\theoremstyle{theorem}
\newtheorem{thm}{Theorem}
\newtheorem{prop}{Proposition}
\newtheorem{lem}{Lemma}
\newtheorem{cor}{Corollary}

\newenvironment{pf}{\begin{proof}}{\end{proof}}
\usepackage{tikz}
\usetikzlibrary{calc,intersections,patterns,through}
\newcommand{\ang}{\tikz@ang}
\def\tikz@ang(#1)(#2)#3{%
\pgfmathanglebetweenpoints{%
\pgfpointanchor{#1}{center}}{%
\pgfpointanchor{#2}{center}}
\pgfmathsetmacro{#3}{\pgfmathresult}%
}

\title{Rational angle bisectors on the coordinate plane\\
and solutions of Pell's equations}
\author{Takashi HIROTSU}
\date{\today}
\begin{document}
\maketitle
\begin{abstract}
On the coordinate plane, the slopes $a$ and $b$ of two straight lines and the slope $c$ of one of their angle bisectors satisfy the equation $(a-c)^2(b^2+1) = (b-c)^2(a^2+1).$ 
Recently, an explicit formula for nontrivial integral solutions of this equation with solutions of negative Pell's equations was discovered by the author. 
In this article, for a given square-free integer $d > 1$ and a given integer $z > 1,$ we describe every integral solution $(x,y)$ of $|x^2-dy^2| = z$ such that $x$ and $dy$ are coprime by using the fundamental unit of $\mathbb Q(\sqrt d)$ and elements of $\mathbb Z[\sqrt d]$ whose absolute value of norms are the smallest prime powers. 
We also describe every nontrivial rational solution of the above equation as one of its applications.
\end{abstract}
\section{Introduction}\label{sec-intro}
On the coordinate plane, we consider the following problem, which can be called
{\itshape the rational angle bisection problem}.
\begin{q}\label{q-bisec}
For which rational numbers $a$ and $b$ are the slopes of the angle bisectors between two straight lines with slopes $a$ and $b$ rational?
\end{q}
\begin{rmk}
Given two straight lines, we consider the two angles formed by them, regardless of whether they are acute or not. 
The bisector of one of the angles and that of the supplementary angle are perpendicular to each other. 
In the case when they are not parallel to the coordinate axes, if one of the slopes is rational, then so is the other, since the product of the slopes is $-1.$
\end{rmk}
Essentially, Problem \ref{q-bisec} has the meaning when the bisectors of $\angle AOB$ can be drawn by connecting $O$ and other lattice points for given lattice points $O,$ $A,$ and $B,$ which is important in drawing techniques. 
Furthermore, in engineering, we can specify the radiation range and the axis of light with a ratio of integers without errors due to approximations to irrational numbers by using a solution to this problem. 
Problem \ref{q-bisec} is reduced to solving the equation 
\begin{equation} 
(a-c)^2(b^2+1) = (b-c)^2(a^2+1) \tag{$\star$} \label{star}  
\end{equation} 
in $\mathbb Q$ (see \cite[Proposition 1]{hir}). 
We say that a solution $(a,b,c)$ of \eqref{star} is {\it trivial} if $|a| = |b|.$
\begin{eg}
The triples 
\[ (a,b,c) = \left(\frac{3}{4},\frac{12}{5},\frac{9}{7}\right),\ \left(\frac{1}{7},\frac{23}{7},\frac{6}{7}\right),\ (1,7,2)\] 
satisfy \eqref{star}, and therefore $(a,b) = (3/4,12/5),$ $(1/7,23/7),$ $(1,7)$ are solutions to Problem \ref{q-bisec} (see the following figures).
\end{eg}
\begin{center}\includegraphics{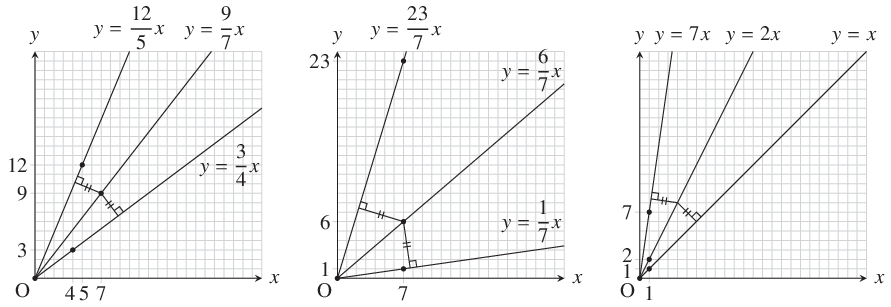}\end{center}\par
The following proposition plays an important role in solving Problem \ref{q-bisec}. \begin{prop}\label{prop-why-pel}
Every nontrivial rational solution $(a,b,c)$ of \eqref{star} is given by 
\begin{align*} 
(a,b,c) = \left( a_1,b_1,\frac{a_1b_2+a_2b_1}{b_2+a_2}\right),\ \left( a_1,b_1,\frac{a_1b_2-a_2b_1}{b_2-a_2}\right) 
\end{align*} 
for some rational solutions $(x,y) = (a_1,a_2),$ $(b_1,b_2)$ of 
\[ x^2-dy^2 = -1,\] 
where $d$ is a positive square-free integer.
\end{prop}
The proof is given in Section 3. 
The following formula for integral solutions of \eqref{star} was proven by using certain criteria for the divisibility of integral solutions of Pell's equations.
\begin{thm}[{\cite[Theorem 1]{hir}}]
For each square-free integer $d > 1$ such that $x^2-dy^2 = -1$ has an integral solution, we denote the $n$-th smallest positive integral solution of $|x^2-dy^2| = 1$ by $(x,y) = (f_n^{(d)},g_n^{(d)}).$ 
Every nontrivial integral solution $(a,b,c)$ of \eqref{star} is given by 
\begin{align} 
(a,b,c) = &\pm\left( f_{(2m-1)(2n-1)}^{(d)},f_{(2m-1)(2n+1)}^{(d)},\frac{g_{(2m-1)\cdot 2n}^{(d)}}{g_{2m-1}^{(d)}}\right), \label{eq-int-d} \\
&\pm (f_{2n-1}^{(2)},-f_{2n+1}^{(2)},f_{2n}^{(2)}) \label{eq-int-2} 
\end{align} 
for some integers $d,$ $m,$ $n > 0,$ after switching $a$ and $b$ if necessary, where \eqref{eq-int-d} contains the case when $d = 2.$ 
Conversely, every triple $(a,b,c)$ of Form \eqref{eq-int-d} or \eqref{eq-int-2} is an integral solution of \eqref{star}. 
\end{thm}
Throughout this article, we use the following definition and symbols.
\begin{df}
Let $d > 1$ be a square-free integer, and let $z > 0$ be an integer.
\begin{enumerate}
\item[(1)]
Let $(x,y)$ be an integral solution of $|x^2-dy^2| = z.$ 
We say that $(x,y)$ is {\itshape strictly primitive}, if $x$ and $dy$ are coprime.
\item[(2)]
Let $(x,y) = (a_1,a_2),$ $(b_1,b_2)$ be positive integral solutions of $|x^2-dy^2| = z.$ 
We say $(a_1,a_2)$ is {\itshape smaller} than $(b_1,b_2),$ if $a_2 < b_2,$ or if $a_1 < b_1$ and $a_2 = b_2.$
\item[(3)]
For each equation of $x^2-dy^2 = z,$ $x^2-dy^2 = -z,$ and $|x^2-dy^2| = z,$ we call its minimum positive integral solution its {\itshape fundamental solution}.
\end{enumerate}
\end{df}
\begin{rmk}
Suppose that $z$ is a square number, and let $(x,y)$ be an integral solution of $|x^2-dy^2| = z.$ 
If $(x,y)$ is {\itshape primitive}, that is, $x$ and $y$ are coprime, then $(x,y)$ is strictly primitive; otherwise, there exists a common prime divisor of $x$ and $d$ dividing $x^2-dy^2 = \pm z$ twice and also $x^2\mp z = dy^2$ twice, which contradicts that $x$ and $y$ are coprime, and $d$ is square-free. 
In particular, every integral solution of $|x^2-dy^2| = 1$ is strictly primitive.
\end{rmk}
\begin{note}\label{note}
Let $\mathbb N$ denote the additive monoid of all nonnegative integers. 
For each prime number $p,$ let $\mathrm{ord}_p:\mathbb Q^\times\to\mathbb Z$ be the normalized $p$-adic additive valuation. 
We denote the greatest common divisor of $a,$ $b \in \mathbb Z\setminus\{ 0\}$ by $\mathrm{gcd}\,(a,b).$
\end{note}
\begin{note}
Let $d > 1$ be a square-free integer. 
Let $\eta$ be the fundamental unit of the real quadratic field $\mathbb Q(\sqrt d).$ 
Let $(x,y) = (f_1,g_1)$ be the fundamental solution of $|x^2-dy^2| = 1,$ and let $\varepsilon = f_1+g_1\sqrt d.$ 
Let 
\[ S(d) = \{ p\text{: prime number} \mid \text{$|x^2-dy^2| = p^l$ has a strictly primitive integral solution for some integer $l > 0$}\}.\] 
For each $p \in S(d),$ let 
\[ l_p = \min\{ l \in \mathbb Z \mid \text{$|x^2-dy^2| = p^l$ has a strictly primitive integral solution},\ l > 0\},\]
and let $(x,y) = (x_p,y_p)$ be the fundamental solution of 
\[\begin{cases} 
x^2-dy^2 = p^{l_p} & \text{if }x^2-dy^2 = -1\text{ has a integral solution,} \\ 
|x^2-dy^2| = p^{l_p} & \text{otherwise}, 
\end{cases}\] 
and let 
\[\xi _p = x_p+y_p\sqrt d.\] 
Let 
\[ S(d)_- = \{ p \in S(d) \mid x_p{}^2-dy_p{}^2 = -p^{l_p}\}.\] 
Furthermore, for each $\alpha = a_1+a_2\sqrt d \in \mathbb Q(\sqrt d),$ where $a_1,$ $a_2 \in \mathbb Q,$ we denote the conjugate and norm of $\alpha$ by 
\[\alpha ' = a_1-a_2\sqrt d \quad\text{and}\quad N(\alpha ) = \alpha\alpha ' = a_1{}^2-da_2{}^2,\] 
respectively.
\end{note}
The first main theorem of this article is the following formula.
\begin{thm}\label{thm-rat-bisec}
Let $(a,b,c)$ be a nontrivial rational solution of \eqref{star}.
\begin{enumerate}
\item[\textup{(i)}]
Suppose that $a$ and $b$ are the $x$-components of rational solutions of the equation $x^2-y^2 = -1.$ 
Then $(a,b,c)$ is given by 
\begin{equation} 
(a,b,c) = \left(\frac{l^2-n^2}{2ln},\frac{m^2-n^2}{2mn},\frac{lm-n^2}{(l+m)n}\right),\ \left(\frac{l^2-n^2}{2ln},\frac{m^2-n^2}{2mn},-\frac{(l+m)n}{lm-n^2}\right) \label{rat-pyt} 
\end{equation} 
for some $l,$ $m,$ $n \in \mathbb Z$ such that $|l| \neq |m|,$ $lm \neq n^2,$ and $lmn \neq 0.$
\item[\textup{(ii)}]
Suppose that $a$ and $b$ are the $x$-components of rational solutions of a common negative Pell's equation $x^2-dy^2 = -1$ for some square-free integer $d > 1.$ 
Then $(a,b,c)$ is given by 
\begin{equation} 
(a,b,c) = \left(\frac{\alpha +\alpha '}{2},\frac{\beta +\beta '}{2},\frac{\alpha\beta -(\alpha\beta )'}{\alpha +\beta -(\alpha  +\beta )'}\right),\ \left(\frac{\alpha +\alpha '}{2},\frac{\beta +\beta '}{2},-\frac{\alpha +\beta -(\alpha +\beta )'}{\alpha\beta -(\alpha\beta )'}\right) \label{rat-pel} 
\end{equation} 
for some $\alpha,$ $\beta \in \mathbb Q(\sqrt d)$ such that $N(\alpha ) = N(\beta ) = -1$ and $\beta \neq \pm\alpha.$ 
In addition, $\alpha$ and $\beta$ can be written in the form 
\[\alpha = \pm\eta ^m\prod_{p \in S(d)}\alpha _p{}^{m_p}p^{-l_pm_p/2} \quad\text{and}\quad \beta = \pm\eta ^n\prod_{p \in S(d)}\beta _p{}^{n_p}p^{-l_pn_p/2}\] 
for some $m_p,$ $n_p \in \mathbb N,$ $\alpha _p,$ $\beta _p \in \{\xi _p,\xi _p'\}$ $(p \in S(d)),$ and $m,$ $n \in \mathbb Z$ satisfying 
\[ l_pm_p \equiv l_pn_p \equiv 0 \pmod 2\] 
and 
\[ 1 \equiv \begin{cases} 
m \equiv n \qquad\qquad\qquad\quad\;\:\,(\mathrm{mod}\ 2) & \text{if }x^2-dy^2 = -1\ \text{has an integral solution}, \\ 
\displaystyle\sum_ {p \in S(d)_-}m_p \equiv \sum_ {p \in S(d)_-}n_p\ (\mathrm{mod}\ 2) & \text{otherwise}. 
\end{cases}\]
\end{enumerate}
Conversely, every triple $(a,b,c)$ of Form \eqref{rat-pyt} or \eqref{rat-pel} is a rational solution of \eqref{star}.
\end{thm}
The proof is given in Section \ref{sec-rat-bisec}. 
\begin{eg}
\begin{enumerate}
\item[(1)]
In \eqref{rat-pyt}, by letting $(l,m,n) = (2,3,1),$ we obtain the rational solutions 
\[ (a,b,c) = \left(\frac{3}{4},\frac{4}{3},1\right),\ \left(\frac{3}{4},\frac{4}{3},-1\right)\] 
of \eqref{star}, and by letting $(l,m,n) = (2,5,1),$ we obtain the rational solutions 
\[ (a,b,c) = \left(\frac{3}{4},\frac{12}{5},\frac{9}{7}\right),\ \left(\frac{3}{4},\frac{12}{5},-\frac{7}{9}\right)\] 
of \eqref{star}.
\item[(2)]
In \eqref{rat-pel}, by letting $d = 2,$ $\alpha = (1+5\sqrt 2)/7,$ $\beta = \eta ^2\alpha,$ and $\eta = 1+\sqrt 2,$ we obtain the rational solutions 
\[ (a,b,c) = \left(\frac{1}{7},\frac{23}{7},\frac{6}{7}\right),\ \left(\frac{1}{7},\frac{23}{7},-\frac{7}{6}\right)\] 
of \eqref{star}, and by letting $d = 34,$ $\alpha = (5+\sqrt{34})/3,$ $\beta = \eta\alpha,$ and $\eta = 35+6\sqrt{34},$ we obtain the rational solutions 
\[ (a,b,c) = \left(\frac{5}{3},\frac{379}{3},\frac{32}{9}\right),\ \left(\frac{5}{3},\frac{379}{3},-\frac{9}{32}\right)\] 
of \eqref{star}.
\end{enumerate}
\end{eg}
The result in (i) refines the previous result \cite[Theorem 2]{hir}. 
For any rational solution $(a,b,c)$ of \eqref{star} in (i), each of $a$ and $b$ is a ratio of the leg lengths of some Pythagorean triangle. 
We find the rational solutions of \eqref{star} in (ii) by using the following Theorem. 
This is the second main theorem of this article.
\begin{thm}\label{thm-gen-pel}
\begin{enumerate}
\item[{\rm (1)}]
Let $z > 1$ be an integer. 
Then $|x^2-dy^2| = z$ has a strictly primitive integral solution if and only if 
\begin{equation} 
\mathrm{ord}_p(z) = \left\{\begin{array}{lll} 
l_pn_p & \text{if }p \in S(d), \\ 
0 & \text{if }p \notin S(d) 
\end{array}\right. \label{eq-cond-idx-spl} 
\end{equation} 
for some $n_p \in \mathbb N$ for each prime number $p,$ where $l_2 = 2$ and $n_2 \in \{ 0,1\}$ if $\eta \notin \mathbb Z[\sqrt d].$ 
In this case, its strictly primitive integral solution $(x,y)$ satisfies 
\begin{equation} 
x+y\sqrt d = \pm\eta ^n\prod_{p \in S(d)}\xi _p^*{}^{n_p} \label{eq-prim-pel} 
\end{equation} 
for some $\xi _p^* \in \{\xi _p,\xi _p'\}$ $(p \in S(d))$ and $n \in \mathbb Z$ such that 
\begin{equation} 
n \equiv \begin{cases} 
0,\ \pm 1 \pmod 3 & \text{if }\eta \in \mathbb Z[\sqrt d]\text{ or \ \;}z \equiv 0\ (\mathrm{mod}\ 2), \\ 
0\qquad\:\!\pmod 3 & \text{if }\eta \notin \mathbb Z[\sqrt d]\text{ and }z \equiv 1\ (\mathrm{mod}\ 2). 
\end{cases} \label{eq-n-trinity} 
\end{equation}
\item[\rm{(2)}]
Let $z > 1$ be an integer. 
Then every integral solution $(x,y)$ of $|x^2-dy^2| = z^2$ satisfies 
\begin{equation} 
x+y\sqrt d = \pm\eta ^n\prod_{p \in S(d)}\xi _p^*{}^{n_p}p^{\mathrm{ord}_p(z)-l_pn_p/2}\prod_{p \notin S(d)}p^{\mathrm{ord}_p(z)} \label{eq-int-pel} 
\end{equation} 
for some $n_p \in \mathbb N,$ $\xi _p^* \in \{\xi _p,\xi _p'\}$ $(p \in S(d)),$ and $n \in \mathbb Z,$ where 
\begin{equation} 
l_pn_p \equiv 0 \pmod 2, \quad n_p \leq 2\,\mathrm{ord}_p(z)/l_p, \label{eq-cond-np} 
\end{equation} 
and \eqref{eq-n-trinity} hold.
\end{enumerate}
\end{thm}
The proof is given in Section \ref{sec-rat-pel}.
\section{Rational Solutions of Pell's Equations}\label{sec-rat-pel}
In this section, we prove Theorem \ref{thm-gen-pel}. 
Throughout the section, let $d > 1$ be a square-free integer. 
We recall the following classical theorems.
\begin{thm}[{\cite[Theorem 3.2.1]{aa}}]
The positive Pell's equation $x^2-dy^2 = 1$ has a nontrivial integral solution independently of the value of $d.$
\end{thm}
\begin{thm}
The negative Pell's equation $x^2-dy^2 = -1$ has a rational solution, if and only if $d$ has no prime divisors congruent to $3$ modulo $4.$
\end{thm}
\begin{pf}
By Fermat's two squares theorem (see \cite[Theorem 2.15]{nzm}), $d$ can be expressed as $d = a^2+b^2$ for some $a,$ $b \in \mathbb N,$ if and only if $d$ has no prime divisors congruent to $3$ modulo $4.$ 
We can assume $b > 0,$ since $d > 1.$ 
Under these conditions, $(x,y) = (a/b,1/b)$ is a rational solution of $x^2-dy^2 = -1.$\par
Conversely, if $x^2-dy^2 = -1$ has a rational solution $(x,y) = (a/c,b/c)$ with $a,$ $b,$ $c \in \mathbb Z$ and $c > 0,$ then $db^2 = a^2+c^2,$ which implies that $db^2$ has no prime divisors congruent to $3$ modulo $4$ by Fermat's theorem as above, and so does $d.$\qedhere
\end{pf}
\begin{rmk}
If $x^2-dy^2 = -1$ has an integral solution, then $d$ has no prime divisors congruent to $3$ modulo $4;$ however, the converse is false. 
For example, $x^2-34y^2 = -1$ has a rational solution $(x,y) = (5/3,1/3)$ but no integral solutions.
\end{rmk}
Let $z,$ $w \in \mathbb Z\setminus\{ 0\}.$ 
If $x^2-dy^2 = z$ and $x^2-dy^2 = w$ have integral solutions $(x,y) = (a_1,a_2)$ and $(x,y) = (b_1,b_2),$ respectively, then $x^2-dy^2 = zw$ has an integral solution 
\[ (x,y) = (a_1b_1+da_2b_2,a_1b_2+a_2b_1),\] 
since the norm map $N:\mathbb Q(\sqrt d)^\times \to \mathbb Q^\times$ is a group homomorphism. 
In the case when $w = 1,$ such a solution is called a {\it Pell multiple}.
\begin{thm}[{\cite[Corollary 3.5]{con}}]
Let $z \in \mathbb Z\setminus\{ 0\}.$ 
Suppose that $x^2-dy^2 = z$ has an integral solution. 
Then, for this equation, there exists a finite number of integral solutions $(x,y) = (a_{1,1},a_{1,2}),$ $\ldots,$ $(a_{r,1},a_{r,2})$ such that every integral solution $(x,y)$ satisfies 
\[ x+y\sqrt d = \pm\varepsilon{}^n(a_{i,1}+a_{i,2}\sqrt d)\] 
for some $i \in \{ 1,\dots,r\}$ and $n \in \mathbb Z.$
\end{thm}
In this article, we generalize the concept of Pell multiples. 
The following propositions are fundamental.
\begin{prop}\label{prop-prod-prim}
\begin{enumerate}
\item[{\rm (1)}]
Let $p_1,$ $\ldots,$ $p_r$ be distinct prime numbers, and let $n_1,$ $\ldots,$ $n_r \in \mathbb N.$ 
If $|x^2-dy^2| = p_i{}^{n_i}$ has a strictly primitive integral solution $(x,y) = (a_{i,1},a_{i,2})$ for each $i \in \{ 1,\dots,r\},$ then $(x,y) \in \mathbb Z^2$ defined by 
\begin{equation} 
x+y\sqrt d = \pm\prod_{i = 1}^r(a_{i,1}+a_{i,2}\sqrt d) \label{eq-prod-prim} 
\end{equation} 
is a strictly primitive integral solution of $|x^2-dy^2| = \prod_{i = 1}^{r}p_i{}^{n_i}.$
\item[{\rm (2)}]
Let $p$ be a prime number, and let $m,$ $n > 0$ be integers. 
If $|x^2-dy^2| = p^m$ has a strictly primitive integral solution $(x,y) = (a_1,a_2),$ then $(x,y) \in \mathbb Z^2$ defined by 
\begin{equation} 
x+y\sqrt d = \pm (a_1+a_2\sqrt d)^n \label{eq-pw-prim} 
\end{equation} 
is a strictly primitive integral solution of $|x^2-dy^2| = p^{mn}.$
\end{enumerate}
\end{prop}
\begin{pf}
We prove the contrapositions.
\begin{enumerate}
\item[(1)]
For each $i \in \{ 1,\dots,r\},$ let $(x,y) = (a_{i,1},a_{i,2})$ be an integral solution of $|x^2-dy^2| = p_i{}^{n_i},$ and let $\alpha _i = a_{i,1}+a_{i,2}\sqrt d.$ 
Suppose that the integral solution $(x,y)$ of $|x^2-dy^2| = p_1{}^{n_1}\cdots p_r{}^{n_r}$ defined by \eqref{eq-prod-prim} is not strictly primitive. 
Then $p_j \mid x$ and $p_j \mid dy$ for some $j \in \{ 1,\dots,r\}.$
\begin{itemize}
\item
{\it Case 1:} Suppose that $p_j \mid d.$ 
Then $p_j \mid da_{j,2}{}^2\pm p_j{}^{n_j} = a_{j,1}{}^2,$ and therefore $p_j \mid a_{j,1}.$
\item
{\it Case 2:} Suppose that $p_j \nmid d.$ 
Then $p_j \mid y.$ 
In $\mathbb Z[\sqrt d],$ this implies that $p_j$ divides 
\[ (x+y\sqrt d)\prod_{i \neq j}\alpha _i' = \pm\alpha _j\prod_{i \neq j}\alpha _i\alpha _i' = \pm\alpha _j\prod_{i \neq j}p_i{}^{n_i}\] 
and therefore $\alpha _j.$ 
Hence $p_j \mid a_{j,1}$ and $p_j \mid a_{j,2}.$
\end{itemize}
Thus, in each case, $p_j \mid a_{j,1}$ and $p_j \mid da_{j,2},$ which imply that the solution $(x,y) = (a_{j,1},a_{j,2})$ is not strictly primitive.
\item[(2)]
Let $(x,y) = (a_1,a_2)$ be an integral solution of $|x^2-dy^2| = p^m.$ 
Suppose that the solution $(x,y)$ of $|x^2-dy^2| = p^{mn}$ defined by \eqref{eq-pw-prim} is not strictly primitive. 
Then $p \mid x$ and $p \mid dy.$
\begin{itemize}
\item
{\it Case 1:} Suppose that $p \mid d.$ 
Then $p \mid da_2{}^2\pm p^m = a_1{}^2,$ and therefore $p \mid a_1.$
\item
{\it Case 2:} Suppose that $p \nmid d.$ 
Then $p \mid y,$ and therefore $p \mid x+y\sqrt d$ in $\mathbb Z[\sqrt d].$ 
With respect to an extension $v$ of $\mathrm{ord}_p:\mathbb Q^\times\to\mathbb Z$ to $\mathbb Q(\sqrt d),$ this implies 
\[ v(a_1+a_2\sqrt d) = \frac{1}{n}v(x+y\sqrt d) > 0.\] 
Hence $p \mid a_1$ and $p \mid a_2.$
\end{itemize}
Thus, in each case, $p \mid a_1$ and $p \mid da_2,$ which imply that the solution $(x,y) = (a_1,a_2)$ is not strictly primitive.\qedhere
\end{enumerate}
\end{pf}
\begin{prop}
Let $z \in \mathbb Z\setminus\{ 0\}.$ 
Suppose that $x^2-dy^2 = -1$ has an integral solution. 
If $x^2-dy^2 = z$ has a strictly primitive integral solution, then so does $x^2-dy^2 = -z.$
\end{prop}
\begin{pf}
We prove the contraposition. 
Suppose that $x^2-dy^2 = -z$ has no strictly primitive integral solutions. 
Let $(x,y) = (x_0,y_0)$ be an integral solution of $x^2-dy^2 = z.$ 
Then $(x,y) \in \mathbb Z^2$ defined by $x+y\sqrt d = \varepsilon (x_0+y_0\sqrt d)$ is an integral solution of $x^2-dy^2 = -z,$ which implies that there exists a prime divisor $p$ of $z$ such that $p \mid x$ and $p \mid dy.$
\begin{itemize}
\item
{\it Case 1:} Suppose that $p \mid d.$ 
Then $p \mid dy_0{}^2+z = x_0{}^2,$ and therefore $p \mid x_0.$
\item
{\it Case 2:} Suppose that $p \nmid d.$ 
Then $p \mid y.$ 
In $\mathbb Z[\sqrt d],$ this implies $p \mid x+y\sqrt d,$ and therefore $p \mid x_0+y_0\sqrt d$ since $p \nmid \varepsilon.$ 
Hence $p \mid x_0$ and $p \mid y_0.$
\end{itemize}
Thus, in each case, $p \mid x_0$ and $p \mid dy_0,$ which imply that $x^2-dy^2 = z$ has no strictly primitive integral solutions.
\end{pf}
We say that an ideal $\mathfrak a \neq (0)$ in $\mathbb Z[\sqrt d]$ or the integer ring $O_K$ of a quadratic field $K$ is {\it primitive} if $\mathfrak a$ is not divisible by $(p)$ for any prime number $p.$
\begin{thm}\label{thm-prim-dec}
Let $p$ be a prime number. 
Suppose that $d \not\equiv 5\ (\mathrm{mod}\ 8)$ or $p \neq 2.$ 
Then $|x^2-dy^2| = z$ has a strictly primitive integral solution for some multiple $z$ of $p,$ if and only if $p$ splits in $\mathbb Q(\sqrt d).$
\end{thm}
\begin{pf}
Let $K = \mathbb Q(\sqrt d).$ 
Recall that the integer ring $O_K$ of $K$ is 
\[ O_K = \begin{cases} 
\mathbb Z\left[\dfrac{1+\sqrt d}{2}\right] & \text{if }d \equiv 1\quad\ \; (\mathrm{mod}\ 4), \\ 
\mathbb Z[\sqrt d] & \text{if }d \equiv 2,\ 3\ (\mathrm{mod}\ 4). 
\end{cases}\] 
Let $I$ be the monoid of all nonzero ideals in $\mathbb Z[\sqrt d].$ 
For each $(\mathfrak a,\mathfrak b) \in I\times I,$ we define $\mathfrak a \sim \mathfrak b$ if 
\[\mathfrak b = (\lambda )\mathfrak a \quad\text{or}\quad \mathfrak a = (\lambda )\mathfrak b\] 
for some $\lambda \in \mathbb Z[\sqrt d]\setminus\{ 0\}.$ 
Then the relation $\sim$ on $I$ is an equivalence relation which is compatible with the multiplication. 
Since the ideal class group $Cl_K$ of $K$ is finite, and so is the quotient monoid $I/\!\sim$ which can be regarded as its submonoid. 
Furthermore, $2$ splits in $K$ if $d \equiv 1\ (\mathrm{mod}\ 8),$ and $2$ is unramified in $K$ if $d \equiv 5\ (\mathrm{mod}\ 8).$\par
Suppose that a prime number $p$ splits into prime ideals $\mathfrak p$ and $\mathfrak p ' = \{\alpha ' \mid \alpha \in \mathfrak p\}$ in $O_K,$ that is, $(p) = \mathfrak p\mathfrak p ' \neq \mathfrak p^2.$ 
Then there exist $l,$ $x,$ $y \in \mathbb Z$ such that 
\begin{equation} 
(\mathfrak p\cap \mathbb Z[\sqrt d])^l = (x+y\sqrt d) \quad\text{and}\quad 0 < l \leq \# Cl_K \label{eq-est-l} 
\end{equation} 
by the finiteness of $I/\!\sim,$ where $x$ and $y$ are coprime since $(\mathfrak p\cap \mathbb Z[\sqrt d])^l$ is primitive. 
Furthermore, $x$ and $d$ are coprime, since $x^2-dy^2 = \pm p^l$ and $p \nmid d.$ 
These imply that \mbox{$|x^2-dy^2| = p^l$} has a strictly primitive integral solution.\par
To prove the converse, suppose that an integer $z > 1$ has a prime divisor $p$ which does not split in $K.$ 
Let $(x,y)$ be an integral solution of $|x^2-dy^2| = z.$
\begin{itemize}
\item
{\it Case 1:} Suppose that $p$ is unramified in $K.$ 
Then $p \neq 2$ or ``$d \equiv 3\ (\mathrm{mod}\ 4)$ and $p = 2$'' by assumption. 
In $O_K,$ the principal ideal $(z)$ is decomposed as 
\[ (z) = (x+y\sqrt d)(x-y\sqrt d),\] 
where both sides are divisible by $(p).$ 
Therefore $(x+y\sqrt d)$ and $(x-y\sqrt d)$ are divisible by $(p)$ in $O_K,$ since the prime ideal $(p)$ divides one of these ideals and also the other because of its self-conjugacy. 
This implies that $(x+y\sqrt d)$ and $(x-y\sqrt d)$ are divisible by $(p)$ in $\mathbb Z[\sqrt d],$ since 
\begin{equation} 
p\,O_K\cap\mathbb Z[\sqrt d] = p\,\mathbb Z[\sqrt d]. \label{eq-pzd} 
\end{equation} 
Hence $p \mid x+y\sqrt d,$ and therefore $p \mid x$ and $p \mid y.$
\item
{\it Case 2:} Suppose that $p$ ramifies in $K.$ 
Then $p \mid d,$ and therefore $p$ divides $dy^2\pm z = x^2,$ which implies $p \mid x.$
\end{itemize}
Thus, in each case, $(x,y)$ is not strictly primitive. 
This proves the desired assertion.
\end{pf}
\begin{cor}
Unless $d \equiv 1\ (\mathrm{mod}\ 4)$ and $p = 2,$ then $l_p$ is not greater than the ideal class number of $\mathbb Q(\sqrt d)$ for each $p \in S(d).$
\end{cor}
\begin{pf}
The desired assertion follows from \eqref{eq-est-l}.
\end{pf}
In the arguments below, we use the following lemma.
\begin{lem}\label{lem}
Suppose that $d \equiv 1\ (\mathrm{mod}\ 4),$ and let $z > 1$ be an integer.
\begin{enumerate}
\item[\textup{(1)}]
If $z \equiv 2\ (\mathrm{mod}\ 4),$ then $|x^2-dy^2| = z$ has no integral solutions. 
Furthermore, if $2 \in S(d),$ then $l_2 \geq 2.$
\item[\textup{(2)}]
Both $2 \in S(d)$ and $l_2 = 2$ hold if and only if $\eta \notin \mathbb Z[\sqrt d].$ 
In this case, $\xi _2 = 2\eta$ holds.
\item[\textup{(3)}]
If $d \equiv 1\ (\mathrm{mod}\ 8)$ and $z \equiv 4\ (\mathrm{mod}\ 8),$ then every integral solution of $|x^2-dy^2| = z$ is a pair of even integers. 
Furthermore, if $d \equiv 1\ (\mathrm{mod}\ 8),$ then $\eta \in \mathbb Z[\sqrt d].$
\item[\textup{(4)}]
If $d \equiv 5\ (\mathrm{mod}\ 8)$ and $z \equiv 0\ (\mathrm{mod}\ 8),$ then every integral solution of $|x^2-dy^2| = z$ is a pair of even integers.
\end{enumerate}
\end{lem}
\begin{pf}
\begin{enumerate}
\item[(1)]
If $z \equiv 2\ (\mathrm{mod}\ 4),$ then $|x^2-dy^2| = z$ has no integral solutions, since $x^2-dy^2 \equiv x^2-y^2 \not\equiv 2\ (\mathrm{mod}\ 4)$ by assumption. 
This implies that $l_2 \geq 2$ if $2 \in S(d).$
\item[(2)]
Suppose that $\eta \notin \mathbb Z[\sqrt d].$ 
Then the pair $(x,y)$ of odd integers defined by $\eta = (x+y\sqrt d)/2$ satisfies $|x^2-dy^2| = 2^2.$ 
This implies that $x$ and $dy$ are coprime, and therefore $2 \in S(d)$ and $l_2 = 2.$ 
In this case, $\xi _2 = 2\eta$ holds by the minimality of $\eta.$\par
Conversely, if $2 \in S(d)$ and $l_2 = 2,$ then a strictly primitive integral solution $(x,y)$ of $|x^2-dy^2| = 2^2$ satisfies $(x+y\sqrt d)/2 \in O_K^\times\setminus\mathbb Z[\sqrt d]^\times,$ which implies $\eta \notin \mathbb Z[\sqrt d]$ since 
\begin{equation} 
O_K^\times = \{\pm\eta ^n \mid n \in \mathbb Z\}. \label{eq-unit-gp} 
\end{equation}
\item[(3)]
Suppose that $d \equiv 1\ (\mathrm{mod}\ 8).$ 
Then 
\[ 1^2-d\cdot 1^2 \equiv 1^2-d\cdot 3^2 \equiv 3^2-d\cdot 1^2 \equiv 3^2-d\cdot 3^2 \equiv 0 \pmod 8,\] 
which imply that every integral solution $(x,y)$ of $|x^2-dy^2| = z$ is a pair of even integers if $z \equiv 4\ (\mathrm{mod}\ 8).$ 
This implies $\eta \in \mathbb Z[\sqrt d]$ by (2).
\item[(4)]
If $d \equiv 5\ (\mathrm{mod}\ 8),$ then 
\[ 1^2-d\cdot 1^2 \equiv 1^2-d\cdot 3^2 \equiv 3^2-d\cdot 1^2 \equiv 3^2-d\cdot 3^2 \equiv 4 \pmod 8,\] 
which imply that every integral solution $(x,y)$ of $|x^2-dy^2| = z$ is a pair of even integers if $z \equiv 0\ (\mathrm{mod}\ 8).$
\end{enumerate}
\end{pf}
Now we can determine the elements of $S(d).$
\begin{thm}\label{thm-sd}
\begin{enumerate}
\item[\textup{(1)}]
If $d \equiv 1\ (\mathrm{mod}\ 8),$ or $d \equiv 5\ (\mathrm{mod}\ 8)$ and $\eta \in \mathbb Z[\sqrt d],$ or $d \equiv 2,$ $3\ (\mathrm{mod}\ 4),$ then $S(d)$ consists of all prime numbers which split in $\mathbb Q(\sqrt d).$
\item[\textup{(2)}]
If $d \equiv 5\ (\mathrm{mod}\ 8)$ and $\eta \notin \mathbb Z[\sqrt d],$ then $S(d)$ consists of all prime numbers which split in $\mathbb Q(\sqrt d)$ and $2.$
\end{enumerate}
\end{thm}
\begin{pf}
If $d \equiv 5\ (\mathrm{mod}\ 8)$ and $\eta \in \mathbb Z[\sqrt d],$ then $2 \notin S(d),$ since $|x^2-dy^2| = 2^l$ has no strictly primitive integral solutions for each integer $l > 0$ by Lemma \ref{lem} (1), (2), and (4). 
Combining this fact with  Theorem \ref{thm-prim-dec} and Lemma \ref{lem} (2), we obtain (1) and (2).\qedhere
\end{pf}
We also use the following fact.
\begin{thm}\label{thm-prim-ideal}
Let $\mathfrak a_1,$ $\dots,$ $\mathfrak a_r$ be pairwise coprime ideals in the integer ring $O_K$ of a quadratic field $K.$ 
Then $\prod _{i = 1}^r\mathfrak a_i$ is primitive if and only if $\mathfrak a_1,$ $\dots,$ $\mathfrak a_r$ are primitive.
\end{thm}
\begin{pf}
For any ideal $\mathfrak a \neq (0)$ in $O_K,$ $\mathfrak a$ is primitive if and only if the residue group $O_K/\mathfrak a$ is cyclic (see \cite[Corollary 6.30]{aok}). 
Combining this fact with the Chinese remainder theorem, we obtain the desired assertion.
\end{pf}
Now we are ready to prove Theorem \ref{thm-gen-pel}.
\begin{pf}[Proof of Theorem \ref{thm-gen-pel}]
\begin{enumerate}
\item[(1)]
Let $K = \mathbb Q(\sqrt d),$ and let $O_K$ denote the integer ring of $K.$ 
For each prime number $p,$ we take a prime ideal $\mathfrak p$ in $O_K$ lying over $(p)$ in $\mathbb Z,$ and let $\mathfrak p{}' = \{\alpha ' \mid \alpha \in \mathfrak p\}.$ 
Note that $\mathfrak p{}^r$ and $\mathfrak p{}'^r$ are primitive for each $r \in \mathbb N$ if $p$ splits in $K.$ 
Let $(x,y)$ be a strictly primitive integral solution of $|x^2-dy^2| = z.$ 
Let $T$ be the set of all prime divisors of $z,$ and let $T(d) = T\cap S(d).$ 
For each $p \in T(d),$ let $q_p$ and $r_p$ be the quotient and remainder, respectively, when dividing $\mathrm{ord}_p(z)$ by $l_p.$ 
For each $p \in T\setminus T(d),$ let $r_p = \mathrm{ord}_p(z).$ 
The equality $|x^2-dy^2| = z$ can be expressed as $(x^2-dy^2) = (z)$ with ideals in $O_K,$ which implies 
\[ (x+y\sqrt d)(x-y\sqrt d) = \prod_{p \in T(d)}(\xi _p)^{q_p}(\xi _p')^{q_p}\prod_{p \in T}(p)^{r_p},\] 
since $(\xi _p)(\xi _p') = (p)^{l_p}$ for each $p \in T(d).$ 
For each $p \in T,$ if $d \equiv 2,$ $3\ (\mathrm{mod}\ 4)$ or $p \neq 2,$ then $(p)$ is decomposed as 
\[ (p) = \mathfrak p\mathfrak p{}' \neq \mathfrak p^2\] 
in $O_K$ by Theorem \ref{thm-prim-dec}, and therefore the ideal $(x+y\sqrt d)$ is divisible by either $\mathfrak p$ or $\mathfrak p{}'$; otherwise, $(p) = \mathfrak p\mathfrak p{}' \mid (x+y\sqrt d)$ in $O_K$ and also $\mathbb Z[\sqrt d]$ by \eqref{eq-pzd}, and therefore $p \mid x$ and $p \mid y,$ which contradict that $x$ and $dy$ are coprime. 
For each $p \in T(d),$ we can see that $(\xi _p)$ is divisible by either $\mathfrak p$ or $\mathfrak p{}'$ by the same argument as above, which implies 
\[\{ (\xi _p),(\xi _p')\} = \begin{cases} 
\{\mathfrak p{}^{l_p},\mathfrak p{}'^{l_p}\} & \text{if }d \equiv 2,\ 3\ (\mathrm{mod}\ 4)\text{ or }p \neq 2, \\ 
\{\mathfrak p{}^{l_2},\mathfrak p{}'^{l_2}\},\ \{ (2)\mathfrak p{}^{l_2-2},(2)\mathfrak p{}'^{l_2-2}\} & \text{if }d \equiv 1\ (\mathrm{mod}\ 4)\text{ and }p = 2. 
\end{cases}\] 
The latter case follows from $(2)^2 = \mathfrak p{}^2\mathfrak p{}'^2 \nmid (\xi _2),$ since 
\begin{equation} 
2^2O_K\cap\mathbb Z[\sqrt d] \subset 2\,\mathbb Z[\sqrt d] \label{eq-4ok-zd} 
\end{equation} 
and the fundamental solution of $|x^2-dy^2| = 2^{l_2}$ is strictly primitive.
\begin{itemize}
\item
{\it Case 1:} Suppose that $d \equiv 2,\ 3\ (\mathrm{mod}\ 4)$ or $2 \notin T.$ 
Then $T = T(d)$ by Theorem \ref{thm-sd} (1), and the ideal $(x+y\sqrt d)$ is decomposed as 
\[ (x+y\sqrt d) = \left(\prod_{p \in T(d)}\xi _p^*{}^{q_p}\right)\prod_{p \in T(d)}\mathfrak p^*{}^{r_p}\] 
with $\xi _p^* \in \{\xi _p,\xi _p'\}$ and $\mathfrak p^* \in \{\mathfrak p,\mathfrak p{}'\}$ in $O_K.$ 
The ideal 
\[\left(\frac{x+y\sqrt d}{\prod_{p \in T(d)}\xi _p^*{}^{q_p}}\right) = \prod_{p \in T(d)}\mathfrak p^*{}^{r_p}\] 
is principal and primitive by Theorem \ref{thm-prim-ideal}, and therefore equal to $O_K,$ which implies that $r_p = 0$ for each $p \in T(d).$ 
Hence there exists $n \in \mathbb Z$ such that 
\[ x+y\sqrt d = \pm\eta ^n\prod_{p \in T(d)}\xi _p^*{}^{q_p}.\] 
Letting $n_p = q_p,$ we obtain \eqref{eq-prim-pel}.
\item
{\it Case 2:} Suppose that $d \equiv 1\ (\mathrm{mod}\ 4)$ and $2 \in T.$\par
Assume that $2 \notin S(d).$ 
Then $T\setminus\{ 2\} = T(d),$ $d \equiv 5\ (\mathrm{mod}\ 8),$ $2$ is unramified in $K,$ $\eta \in \mathbb Z[\sqrt d]$ by Theorem \ref{thm-sd}, and $\mathrm{ord}_2(z) = 2$ by Lemma \ref{lem} (1) and (4). 
The prime ideal $(2)$ divides both $(x+y\sqrt d)$ and $(x-y\sqrt d)$ only once by its self-conjugacy; otherwise, $(2)^2 \mid (x+y\sqrt d)$ in $O_K,$ and therefore $(2) \mid (x+y\sqrt d)$ in $\mathbb Z[\sqrt d]$ by \eqref{eq-4ok-zd}, which contradicts that $x$ and $dy$ are coprime. 
The ideal $((x+y\sqrt d)/2)$ is decomposed as 
\[\left(\frac{x+y\sqrt d}{2}\right) = \left(\prod_{p \in T(d)}\xi _p^*{}^{q_p}\right)\prod_{p \in T(d)}\mathfrak p^*{}^{r_p}\] 
with $\xi _p^* \in \{\xi _p,\xi _p'\}$ and $\mathfrak p^* \in \{\mathfrak p,\mathfrak p{}'\}$ in $O_K.$ 
The ideal 
\[\left(\frac{x+y\sqrt d}{2\prod_{p \in T(d)}\xi _p^*{}^{q_p}}\right) = \prod_{p \in T(d)}\mathfrak p^*{}^{r_p}\] 
is principal and primitive by Theorem \ref{thm-prim-ideal}, and therefore equal to $O_K,$ which implies that $r_p = 0$ for each $p \in T(d).$ 
Hence there exists $n \in \mathbb Z$ such that 
\[ x+y\sqrt d = \pm 2\eta ^n\prod_{p \in T(d)}\xi _p^*{}^{q_p}.\] 
This implies $x+y\sqrt d \in 2\mathbb Z[\sqrt d],$ and therefore $2 \mid x$ and $2 \mid y,$ which contradict that $x$ and $dy$ are coprime.\par
Thus $2 \in S(d),$ and therefore $T = T(d)$ by Theorem \ref{thm-sd}.
\begin{itemize}
\item
{\it Subcase 2-1:} Suppose that $\eta \in \mathbb Z[\sqrt d].$ 
Then $d \equiv 1\ (\mathrm{mod}\ 8),$ $2$ splits in $K$ by Theorem \ref{thm-sd} (1), and $l_2 \geq 3$ by Lemma \ref{lem} (1) and (2). 
By the same argument as \mbox{Case 1}, we can see that $r_p = 0$ for each $p \in T(d).$ 
Hence there exists $n \in \mathbb Z$ such that 
\[ x+y\sqrt d = \pm\eta ^n\prod_{p \in T(d)}\xi _p^*{}^{q_p}.\] 
Letting $n_p = q_p,$ we obtain \eqref{eq-prim-pel}.
\item
{\it Subcase 2-2:} Suppose that $\eta \notin \mathbb Z[\sqrt d].$ 
Then $d \equiv 5\ (\mathrm{mod}\ 8),$ $2$ is unramified in $K$ by Lemma \ref{lem} (3), and $\mathrm{ord}_2(z) = 2,$ $l_2 = 2,$ $q_2 = 1,$ $r_2 = 0,$ and $\xi _2 = 2\eta$ by Lemma \ref{lem} (2) and (4). 
The ideal $((x+y\sqrt d)/2)$ is decomposed as 
\[\left(\frac{x+y\sqrt d}{2}\right) = \left(\prod_{p \in T(d)\setminus\{ 2\}}\xi _p^*{}^{q_p}\right)\prod_{p \in T(d)\setminus\{ 2\}}\mathfrak p^*{}^{r_p}\] 
with $\xi _p^* \in \{\xi _p,\xi _p'\}$ and $\mathfrak p^* \in \{\mathfrak p,\mathfrak p{}'\}$ in $O_K.$ 
The ideal 
\[\left(\frac{x+y\sqrt d}{2\prod_{p \in T(d)\setminus\{ 2\}}\xi _p^*{}^{q_p}}\right) = \prod_{p \in T(d)\setminus\{ 2\}}\mathfrak p^*{}^{r_p}\] 
is principal and primitive by Theorem \ref{thm-prim-ideal}, and therefore equal to $O_K,$ which implies that $r_p = 0$ for each $p \in T(d)\setminus\{ 2\}.$ 
Hence there exists $m \in \mathbb Z$ such that 
\begin{align*} 
x+y\sqrt d &= \pm 2\eta ^m\prod_{p \in T(d)\setminus\{ 2\}}\xi _p^*{}^{q_p} \\ 
&= \pm\eta ^{m-1}\xi _2\prod_{p \in T(d)\setminus\{ 2\}}\xi _p^*{}^{q_p} = \pm\eta ^{m+1}\xi _2'\prod_{p \in T(d)\setminus\{ 2\}}\xi _p^*{}^{q_p}. 
\end{align*} 
Letting $n_2 = 1,$ $n_p = q_p$ for each $p \neq 2,$ and $n = m\pm 1,$ we obtain \eqref{eq-prim-pel}.
\end{itemize}
\end{itemize}
In these cases, the exponent $n$ satisfies \eqref{eq-n-trinity}, since 
\[\mathbb Z[\sqrt d] ^\times = \begin{cases} 
\{\pm\eta ^n \mid n \in \mathbb Z\} & \text{if }\eta \in \mathbb Z[\sqrt d], \\ 
\{\pm\eta ^n \mid n \in \mathbb Z,\ n \equiv 0\ (\mathrm{mod}\ 3)\} & \text{if }\eta \notin \mathbb Z[\sqrt d] 
\end{cases}\] 
by \eqref{eq-unit-gp} and $\eta ^3 \in \mathbb Z[\sqrt d]$ (see \cite[Theorem 2.1.4]{mol}). 
We can verify that \eqref{eq-cond-idx-spl} is satisfied by the argument above. 
Thus the assertion of (1) holds.
\item[(2)]
Let $(x,y)$ be an integral solution of $|x^2-dy^2| = z^2,$ and let $g = \mathrm{gcd}\,(x,y).$ 
Let $T$ be the set of all prime divisors of $zg^{-1},$ and let $T(d) = T\cap S(d)$. 
Let 
\[ x_0 = \frac{x}{g} \quad\text{and}\quad y_0 = \frac{y}{g}.\] 
Then $(x_0,y_0)$ is a strictly primitive integral solution of $|x^2-dy^2| = z^2g^{-2},$ and satisfies 
\begin{equation} 
x_0+y_0\sqrt d = \pm\eta ^n\prod_{p \in S(d)}\xi _p^*{}^{n_p} \label{eq-prim-sq} 
\end{equation} 
for some $n_p \in \mathbb N,$ $\xi _p^* \in \{\xi _p,\xi _p'\}$ $(p \in S(d)),$ and $n \in \mathbb Z$ satisfying \eqref{eq-n-trinity} and the condition obtained by replacing $z$ with $z^2g^{-2}$ in \eqref{eq-cond-idx-spl} because of (1). 
If $p \in S(d),$ then \eqref{eq-cond-np} holds, since $0 \leq \mathrm{ord}_p(g) = \mathrm{ord}_p(z)-l_pn_p/2.$ 
If $p \notin S(d),$ then $\mathrm{ord}_p(g) = \mathrm{ord}_p(z).$ 
Multiplying both sides of \eqref{eq-prim-sq} by 
\[ g = \prod_{p \in T(d)}p^{\mathrm{ord}_p(z)-l_pn_p/2}\prod_{p \in T\setminus T(d)}p^{\mathrm{ord}_p(z)},\] 
we obtain \eqref{eq-int-pel}.\qedhere
\end{enumerate}
\end{pf}
\begin{rmk}
By Theorem \ref{thm-gen-pel}, it is impossible for both $x^2-dy^2 = p^n$ and $x^2-dy^2 = -p^n$ to have strictly primitive integral solutions, if $x^2-dy^2 = -1$ has no integral solutions.
\end{rmk}
\begin{eg}
The equation $x^2-34y^2 = 1$ has the fundamental solution $(x,y) = (35,6);$ however, $x^2-34y^2 = -1$ has no integral solutions. 
The fundamental unit of $\mathbb Q(\sqrt{34})$ is $\eta = 35+6\sqrt{34}.$ 
The equations $x^2-34y^2 = -3^2,$ $x^2-34y^2 = -5^2,$ and $x^2-34y^2 = -11^2$ have the fundamental solutions $(x,y) = (5,1),$ $(x,y) = (3,1),$ and $(x,y) = (27,5),$ respectively. 
Every integral solution $(x,y)$ of $x^2-34y^2 = -(3\cdot 5\cdot 11)^2$ satisfies one of the following conditions for some $n \in \mathbb Z.$
\begin{itemize}
\item
$x+y\sqrt d = \pm\eta ^n\cdot (5\pm\sqrt{34})\cdot 5\cdot 11.$
\item
$x+y\sqrt d = \pm\eta ^n\cdot 3\cdot (3\pm\sqrt{34})\cdot 11.$
\item
$x+y\sqrt d = \pm\eta ^n\cdot 3\cdot 5\cdot (27\pm5\sqrt{34}).$
\item
$x+y\sqrt d = \pm\eta ^n\cdot (5\pm\sqrt{34})\cdot (3\pm\sqrt{34})\cdot (27\pm5\sqrt{34}).$
\end{itemize}
\end{eg}
Theorem \ref{thm-gen-pel} implies the following formula for rational solutions of Pell's equations.
\begin{thm}\label{thm-rat-pel}
Let $r \in \{ 0,1\}.$ 
Every rational solution $(x,y)$ of $x^2-dy^2 = (-1)^r$ satisfies 
\begin{equation} 
x+y\sqrt d = \pm\eta ^n\prod_{p \in S(d)}\xi _p^*{}^{n_p}p^{-l_pn_p/2} \label{eq-rat-pel} 
\end{equation} 
for some $n_p \in \mathbb N,$ $\xi _p^* \in \{\xi _p,\xi _p'\}$ $(p \in S(d)),$ and $n \in \mathbb Z$ such that 
\[ l_pn_p \equiv 0 \pmod 2\] 
and 
\begin{equation} 
r \equiv \begin{cases} 
n\qquad\quad\:\:\:\pmod 2 & \text{if }x^2-dy^2 = -1\ \text{has an integral solution}, \\ 
\displaystyle\sum_ {p \in S(d)_-}n_p \pmod 2 & \text{otherwise}. 
\end{cases} \label{eq-cond-idx} 
\end{equation} 
\end{thm}
\begin{pf}
Every integral solution $(x,y)$ of $x^2-dy^2 = (-1)^r$ satisfies 
\[ x+y\sqrt d = \pm\eta ^n\] 
for some $n \in \mathbb Z.$\par
Let $(x,y) = (X/Z,Y/Z)$ be a rational solution of $x^2-dy^2 = (-1)^r,$ where $X,$ $Y,$ and $Z$ are coprime integers and $Z > 1.$ 
Then $(X,Y)$ is an integral solution of $X^2-dY^2 = (-1)^rZ^2,$ and satisfies 
\[ X+Y\sqrt d = \pm\eta ^n\prod_{p \in S(d)}\xi _p^*{}^{n_p}p^{\mathrm{ord}_p(Z)-l_pn_p/2}\prod_{p \notin S(d)}p^{\mathrm{ord}_p(Z)}\] 
for some $n_p \in \mathbb N,$ $\xi _p^* \in \{\xi _p,\xi _p'\}$ $(p \in S(d)),$ and $n \in \mathbb Z,$ where 
\[ l_pn_p \equiv 0 \pmod 2, \quad n_p \leq 2\,\mathrm{ord}_p(Z)/l_p,\] 
and 
\begin{equation} 
n \equiv \begin{cases} 
0,\ \pm 1 \pmod 3 & \text{if }\eta \in \mathbb Z[\sqrt d]\text{ or \ \:\,}Z \equiv 0\ (\mathrm{mod}\ 2), \\ 
0\qquad\:\!\pmod 3 & \text{if }\eta \notin \mathbb Z[\sqrt d]\text{ and }Z \equiv 1\ (\mathrm{mod}\ 2) 
\end{cases} \label{eq-add-cond} 
\end{equation} 
by Theorem \ref{thm-gen-pel}. 
Dividing both sides by $Z,$ we obtain \eqref{eq-rat-pel}. 
Condition \eqref{eq-cond-idx} follows from 
\[\begin{cases} 
N(\xi _p^*) < 0 & \text{if }p \in S(d)_-, \\ 
N(\xi _p^*) > 0 & \text{if }p \in S(d)\setminus S(d)_- 
\end{cases}\] 
and 
\[\quad \begin{cases} 
N(\eta ) < 0 & \text{if }x^2-dy^2 = -1\ \text{has an integral solution}, \\ 
N(\eta ) > 0 & \text{otherwise}. 
\end{cases}\] 
We can remove Condition \eqref{eq-add-cond}, since $x,$ $y \in \mathbb Q$ and $\eta ^n \in \mathbb Z[(1+\sqrt d)/2]$ for each $n \in \mathbb Z.$ 
Thus the assertion of the theorem holds.
\end{pf}
\section{Rational Angle Bisectors}\label{sec-rat-bisec}
In this section, we prove Proposition \ref{prop-why-pel} and Theorem \ref{thm-rat-bisec}.
\begin{pf}[{\bfseries Proof of Proposition \ref{prop-why-pel}}]
Let $(a,b,c) = (a_1,b_1,c_1)$ be a nontrivial rational solution of \eqref{star}, and let 
\[ a_1 = \dfrac{A_1}{Z}, \quad b_1 = \dfrac{B_1}{Z}, \quad\text{and}\quad c_1 = \dfrac{C_1}{Z}\] 
with $A_1,$ $B_1,$ $C_1,$ $Z \in \mathbb Z$ and $Z \neq 0.$ 
Substituting these into \eqref{star} and multiplying both sides by $Z^4,$ we obtain  
\begin{equation} 
(A_1-C_1)^2(B_1{}^2+Z^2) = (B_1-C_1)^2(A_1{}^2+Z^2). \label{eq-star-homog} 
\end{equation} 
For any prime number $p,$ the parities of $\mathrm{ord}_p(A_1{}^2+Z^2)$ and $\mathrm{ord}_p(B_1{}^2+Z^2)$ coincide with each other, since $(A_1-C_1)^2$ and $(B_1-C_1)^2$ are square numbers. 
Let $d$ be the product of every prime number $p$ such that these valuations are odd (define $d = 1$ if there are no such prime numbers). 
Then there exist $A_2,$ $B_2 \in \mathbb Z$ such that 
\begin{equation} 
A_1{}^2+Z^2 = dA_2{}^2 \quad\text{and}\quad B_1{}^2+Z^2 = dB_2{}^2, \label{eq-pels} 
\end{equation} 
where $A_1$ and $B_1$ are the $X$-components of the integral solutions $(X,Y) = (A_1,A_2),$ $(B_1,B_2)$ of $X^2-dY^2 = -Z^2.$ 
Furthermore, $a_1$ and $b_1$ are the $x$-components of the rational solutions $(x,y) = (a_1,a_2),$ $(b_1,b_2)$ of $x^2-dy^2 = -1$ with $a_2 = A_2/Z$ and $b_2 = B_2/Z.$ 
Substituting \eqref{eq-pels} into \eqref{eq-star-homog} and dividing both sides by $d,$ we obtain 
\[ (A_1-C_1)^2B_2{}^2 = (B_1-C_1)^2A_2{}^2,\] 
or equivalently, 
\[ (A_1-C_1)B_2 = \pm (B_1-C_1)A_2.\] 
Solving for $C_1,$ we obtain 
\[ C_1 = \frac{A_1B_2+A_2B_1}{B_2+A_2}, \frac{A_1B_2-A_2B_1}{B_2-A_2},\] 
since $B_2{}^2-A_2{}^2 = (B_1{}^2-A_1{}^2)/d = (b_1{}^2-a_1{}^2)Z^2/d \neq 0,$ and therefore 
\[ c_1 = \frac{a_1b_2+a_2b_1}{b_2+a_2}, \frac{a_1b_2-a_2b_1}{b_2-a_2}.\qedhere\] 
\end{pf}
\begin{pf}[{\bfseries Proof of Theorem \ref{thm-rat-bisec}}]
\begin{enumerate}
\item[(i)]
Suppose that $a$ and $b$ are the $x$-components of rational solutions $(x,y) = (a_1,a_2)$ and $(x,y) = (b_1,b_2)$ of $x^2-y^2 = -1,$ respectively. 
Let 
\[ (a_1,a_2) = \left(\frac{A_1}{n},\frac{A_2}{n}\right) \quad\text{and}\quad (b_1,b_2) = \left(\frac{B_1}{n},\frac{B_2}{n}\right)\] 
with $A_1,$ $A_2,$ $B_1,$ $B_2,$ $n \in \mathbb Z$ and $n \neq 0.$ 
Then we have 
\[ A_2{}^2-A_1{}^2 = B_2{}^2-B_1{}^2 = n^2.\] 
Furthermore, let 
\[ A_2+A_1 = l \quad\text{and}\quad B_2+B_1 = m.\] 
Then we have 
\[ A_2-A_1 = \frac{n^2}{l} \quad\text{and}\quad B_2-B_1 = \frac{n^2}{m}.\] 
These imply 
\begin{align*} 
a_1 &= \frac{l-n^2/l}{2n} = \frac{l^2-n^2}{2ln}, & a_2 &= \frac{l+n^2/l}{2n} = \frac{l^2+n^2}{2ln}, \\ 
b_1 &= \frac{m-n^2/m}{2n} = \frac{m^2-n^2}{2mn}, & b_2 &= \frac{m+n^2/m}{2n} = \frac{m^2+n^2}{2mn}, 
\end{align*} 
and therefore 
\[\frac{a_1b_2+a_2b_1}{b_2+a_2} = \frac{1}{2}\cdot\frac{(l^2-n^2)(m^2+n^2)+(l^2+n^2)(m^2-n^2)}{ln(m^2+n^2)+mn(l^2+n^2)} = \frac{1}{2}\cdot\frac{2(l^2m^2-n^4)}{n(lm+n^2)(l+m)} = \frac{lm-n^2}{(l+m)n},\] 
which implies 
\[\frac{a_1b_2-a_2b_1}{b_2-a_2} = -\left(\frac{a_1b_2+a_2b_1}{b_2+a_2}\right) ^{-1} = -\frac{(l+m)n}{lm-n^2}.\]
\item[(ii)]
Suppose that $a$ and $b$ are the $x$-components of rational solutions $(x,y) = (a_1,a_2)$ and $(x,y) = (b_1,b_2)$ of $x^2-dy^2 = -1$ for some square-free integer $d > 1,$ respectively. 
In general, if $\alpha = a_1+a_2\sqrt d$ and $\beta = b_1+b_2\sqrt d$ with $a_1,$ $a_2,$ $b_1,$ $b_2 \in \mathbb Z,$ then 
\[ a_1 = \frac{\alpha +\alpha '}{2} \quad\text{and}\quad b_1 = \frac{\beta +\beta '}{2},\] 
which imply 
\[ a_1b_2+a_2b_1 = \frac{\alpha\beta -(\alpha\beta )'}{2\sqrt d} \quad\text{and}\quad a_2+b_2 = \frac{(\alpha +\beta )-(\alpha +\beta )'}{2\sqrt d}.\] 
Combining these identities with Theorem \ref{thm-rat-pel}, we obtain the desired assertion.\qedhere
\end{enumerate}
\end{pf}

\newpage
\section*{Appendix}\label{sec-app}
In Theorem \ref{thm-rat-bisec}, for $d \leq 34$ such that $x^2-dy^2 = -1$ has a rational solution and $p \leq 97,$ the values of $\eta$ and $\xi _p$ are summarized in the following table. 
For reference, we add the values of $N(\eta )$ and $N(\xi _p)$ below the values of $\eta$ and $\xi _p,$ respectively. 
We also add the value of the ideal class number $h$ of $\mathbb Q(\sqrt d).$
{\fontsize{9pt}{9pt}\selectfont\[\begin{array}{|c||c|c|c|c|c|c|c|c|} \hline 
d & 2 & 5 & 10 & 13 & 17 & 26 & 29 & 34 \\ \hline\hline 
h & 1 & 1 & 2 & 1 & 1 & 2 & 1 & 2 \\ \hline 
\eta & 1+\sqrt 2 & (1+\sqrt 5)/2 & 3+\sqrt{10} & (3+\sqrt{13})/2 & 4+\sqrt{17} & 5+\sqrt{26} & (5+\sqrt{29})/2 & 35+6\sqrt{34} \\ 
{} & -1 & -1 & -1 & -1 & -1 & -1 & -1 & 1 \\ \hline 
\xi _2 & {} & 3+\sqrt 5 & {} & 11+3\sqrt{13} & 5+\sqrt{17} & {} & 27+5\sqrt{29} & {} \\ 
{} & {} & 2^2 & {} & 2^2 & 2^3 & {} & 2^2 & {} \\ \hline 
\xi _3 & {} & {} & 7+2\sqrt{10} & 4+\sqrt{13} & {} & {} & {} & 5+\sqrt{34} \\ 
{} & {} & {} & 3^2 & 3 & {} & {} & {} & -3^2 \\ \hline 
\xi _5 & {} & {} & {} & {} & {} & 21+4\sqrt{26} & 11+2\sqrt{29} & 3+\sqrt{34} \\ 
{} & {} & {} & {} & {} & {} & 5^2 & 5 & -5^2 \\ \hline 
\xi _7 & 3+\sqrt 2 & {} & {} & {} & {} & {} & 6+\sqrt{29} & {} \\ 
{} & 7 & {} & {} & {} & {} & {} & 7 & {} \\ \hline 
\xi _{11} & {} & 4+\sqrt 5 & {} & {} & {} & 15+2\sqrt{26} & {} & 27+5\sqrt{34} \\ 
{} & {} & 11 & {} & {} & {} & 11^2 & {} & -11^2 \\ \hline 
\xi _{13} & {} & {} & 23+6\sqrt{10} & {} & 9+2\sqrt{17} & {} & 97+18\sqrt{29} & {} \\ 
{} & {} & {} & 13^2 & {} & 13 & {} & 13 & {} \\ \hline 
\xi _{17} & 5+2\sqrt 2 & {} & {} & 15+4\sqrt{13} & {} & 11+2\sqrt{26} & {} & {} \\ 
{} & 17 & {} & {} & 17 & {} & 17 & {} & {} \\ \hline 
\xi _{19} & {} & 8+3\sqrt 5 & {} & {} & 6+\sqrt{17} & 45+8\sqrt{26} & {} & {} \\ 
{} & {} & 19 & {} & {} & 19 & 19^2 & {} & {} \\ \hline 
\xi _{23} & 5+\sqrt 2 & {} & {} & 6+\sqrt{13} & {} & 7+\sqrt{26} & 38+7\sqrt{29} & {} \\ 
{} & 23 & {} & {} & 23 & {} & 23 & 23 & {} \\ \hline 
\xi _{29} & {} & 7+2\sqrt 5 & {} & 9+2\sqrt{13} & {} & {} & {} & 3+5\sqrt{34} \\ 
{} & {} & 29 & {} & 29 & {} & {} & {} & -29^2 \\ \hline 
\xi _{31} & 7+3\sqrt 2 & 6+\sqrt 5 & 11+3\sqrt{10} & {} & {} & {} & {} & {} \\ 
{} & 31 & 31 & 31 & {} & {} & {} & {} & {} \\ \hline 
\xi _{37} & {} & {} & 53+12\sqrt{10} & {} & {} & 63+10\sqrt{26} & {} & 141+25\sqrt{34} \\ 
{} & {} & {} & 37^2 & {} & {} & 37^2 & {} & -37^2 \\ \hline 
\xi _{41} & 7+2\sqrt 2 & 11+4\sqrt 5 & 9+2\sqrt{10} & {} & {} & {} & {} & {} \\ 
{} & 41 & 41 & 41 & {} & {} & {} & {} & {} \\ \hline 
\xi _{43} & {} & {} & 47+6\sqrt{10} & 76+21\sqrt{13} & 14+3\sqrt{17} & {} & {} & {} \\ 
{} & {} & {} & 43^2 & 43 & 43 & {} & {} & {} \\ \hline 
\xi _{47} & 7+\sqrt 2 & {} & {} & {} & 8+\sqrt{17} & {} & {} & 9+\sqrt{34} \\ 
{} & 47 & {} & {} & {} & 47 & {} & {} & 47 \\ \hline 
\xi _{53} & {} & {} & 143+42\sqrt{10} & 51+14\sqrt{13} & 11+2\sqrt{17} & {} & 13+2\sqrt{29} & {} \\ 
{} & {} & {} & 53^2 & 53 & 53 & {} & 53 & {} \\ \hline 
\xi _{59} & {} & 8+\sqrt 5 & {} & {} & 22+5\sqrt{17} & 85+12\sqrt{26} & 28+5\sqrt{29} & {} \\ 
{} & {} & 59 & {} & {} & 59 & 59^2 & 59 & {} \\ \hline 
\xi _{61} & {} & 9+2\sqrt 5 & {} & 23+6\sqrt{13} & {} & {} & {} & 45+13\sqrt{34} \\ 
{} & {} & 61 & {} & 61 & {} & {} & {} & -61^2 \\ \hline 
\xi _{67} & {} & {} & 77+12\sqrt{10} & {} & 30+7\sqrt{17} & 167+30\sqrt{26} & 264+49\sqrt{29} & {} \\ 
{} & {} & {} & 67^2 & {} & 67 & 67^2 & 67 & {} \\ \hline 
\xi _{71} & 11+5\sqrt 2 & 14+5\sqrt 5 & 9+\sqrt{10} & {} & {} & {} & 10+\sqrt{29} & {} \\ 
{} & 71 & 71 & 71 & {} & {} & {} & 71 & {} \\ \hline 
\xi _{73} & 9+2\sqrt 2 & {} & {} & {} & {} & {} & {} & {} \\ 
{} & 73 & {} & {} & {} & {} & {} & {} & {} \\ \hline 
\xi _{79} & 9+\sqrt 2 & 18+7\sqrt 5 & 13+3\sqrt{10} & 14+3\sqrt{13} & {} & 27+5\sqrt{26} & {} & {} \\ 
{} & 79 & 79 & 79 & 79 & {} & 79 & {} & {} \\ \hline 
\xi _{83} & {} & {} & 173+48\sqrt{10} & {} & 10+\sqrt{17} & 317+60\sqrt{26} & 92+17\sqrt{29} & {} \\ 
{} & {} & {} & 83^2 & {} & 83 & 83^2 & 83 & {} \\ \hline 
\xi _{89} & 11+4\sqrt 2 & 13+4\sqrt 5 & 27+8\sqrt{10} & {} & 19+4\sqrt{17} & {} & {} & 15+2\sqrt{34} \\ 
{} & 89 & 89 & 89 & {} & 89 & {} & {} & 89 \\ \hline 
\xi _{97} & 13+6\sqrt 2 & {} & {} & {} & {} & {} & {} & {} \\ 
{} & 97 & {} & {} & {} & {} & {} & {} & {} \\ \hline 
\end{array}\]}

\begin{thebibliography}{99}
\bibitem{aa}
T.~Andreescu and D.~Andrica, {\it Quadratic Diophantine Equations} ({\it Developments in Mathematics} {\bf 40}), Springer, New York, 2015.
\bibitem{aok}
N.~Aoki, {\it Number Theory of Prime Numbers and Quadratic Fields}, Kyoritsu Publ., Tokyo, 2012, in Japanese.
\bibitem{car}
R.~D.~Carmichael, On the numerical factors of the arithmetic forms $\alpha ^n\pm\beta ^n,$ {\it Ann.~of Math.~(2)}, {\bf 15} (1913--1914), no.~1/4, 30--48.
\bibitem{con}
K.~Conrad, Pell's equation, II, {\tt https://api.semanticscholar.org/CorpusID:14314437} (Retrieved November 15, 2023).
\bibitem{hir}
T.~Hirotsu, Diophantine equation related to angle bisectors and solutions of Pell's equations, {\tt https://arxiv.org/abs/2209.10434}
\bibitem{kos}
T.~Koshy, {\it Pell and Pell--Lucas Numbers with Applications}, Springer, New York, 2014.
\bibitem{mol}
R.~A.~Mollin, {\it Quadratics}, CRC Press, Boca Raton, FL, 1996.
\bibitem{nzm}
I.~Niven, H.~S.~Zuckerman and H.~L.~Montgomery, {\it An Introduction to
the Theory of Numbers}, 5th ed., Wiley, New York, 1991.
\end{thebibliography}
\end{document}